\documentclass[11pt]{amsart}

\newtheorem{theorem}{Theorem}[section]
\newtheorem{lemma}[theorem]{Lemma}

\theoremstyle{definition}

\newtheorem{corollary}[theorem]{Corollary}

\theoremstyle{remark}

\numberwithin{equation}{section}

\def\phi{\varphi}
\def\epsilon{\varepsilon}
\def\RE{\mbox{\rm Re}\,}
\def\IM{\mbox{\rm Im}\,}

\def\CC{{\mathcal C }}

\def\CT{{\mathcal T}}

\def\C{{\mathbb C }}

\begin{document}

\title[Generating Functions for Hausdorff Moment Sequences]
      {A Note on Generating Functions for Hausdorff Moment Sequences}

\author[O. Roth]{Oliver Roth}
\address{Mathematisches Institut,
Universit\"at W\"urzburg\\ D-97074 W\"urzburg,
Germany}
\email{roth@mathematik.uni-wuerzburg.de}

\author[S. Ruscheweyh]{Stephan Ruscheweyh}
\address{Mathematisches Institut,
Universit\"at W\"urzburg\\ D-97074 W\"urzburg,
Germany}
\email{ruscheweyh@mathematik.uni-wuerzburg.de}

\author[L. Salinas]{Luis Salinas}
\address{Departamento de Inform\'atica,
Universidad T\'ecnica F. Santa Mar\'\i a,
Valpara\'\i so, Chile}
\email{lsalinas@inf.utfsm.cl}


\date{}

\commby{Andreas Seeger}

\thanks{ O.R. and S.R. acknowledge partial
  support also from the German-Israeli 
  Foundation  (grant G-809-234.6/2003). S.R. and L.S. received partial
  support from FONDECYT  
(grants 1070269
  and 7070131) and DGIP-UTFSM (grant 240721).}
\date{July 17, 2007}

\begin{abstract}
For functions $f$ whose Taylor coefficients at the origin form a
Hausdorff moment sequence we study the behaviour of $w(y):=|f(\gamma+iy)|$
for $y>0$ ($\gamma\leq1$ fixed). 
\end{abstract}

\maketitle

\section{Introduction and statement of the results}
A sequence $\{a_k\}_{k\geq0}$ of non-negative real numbers, $a_0=1$,
is called a Hausdorff moment sequence if 
there is a
probability measure\footnote{Here and in the sequel we always assume
  that the measures are Borel} $\mu$ on $[0,1]$ such that 
$$
a_k=\int_{0}^{1}t^k\,d\mu(t),\quad k\geq0,
$$
or, equivalently,
$$
F(z)=\sum_{k=0}^{\infty}a_k\,z^k\,=\,\int_{0}^{1}\frac{d\mu(t)}{1-tz},
$$
and $F$ is its generating
function.

It is well known  (Hausdorff \cite{haus}) that a sequence $\{a_k\}_{k\geq0}$ with $a_0=1$ is a Hausdorff moment sequence if and
only if it is  {\em completely monotone }  i.e. 
\begin{eqnarray*}
\Delta^n a_k &:=&\Delta^{n-1}a_k - \Delta^{n-1} a_{k+1}\,\geq\,0,\quad
k\geq0,\ n\geq1,
\end{eqnarray*}
where $\Delta^0$ is the identity operator: $\Delta^0a=a$.

Let $\CT$ denote the set of such generating functions $F$. They are analytic 
in the slit domain  $\Lambda:=\C\setminus[1,\infty)$ and also belong to the
set of Pick functions $P(-\infty,1)$ (see Donoghue \cite{dono} for
more information  on Pick functions).

Wirths \cite{wi} has shown that $f\in\CT$ implies that the function
$zf(z)$ is univalent in
the half-plane $\RE z<1$, and recently the theory of universally
prestarlike mappings has been developed, showing a close link to
$\CT$, see \cite{rss}.  Many classical functions belong to $\CT$ or
are closely related to it. We mention only the polylogarithms
$$
Li_\alpha(z):=\sum_{k=1}^{\infty}\frac{z^k}{k^\alpha},\quad\alpha\geq0,
$$ 
where $Li_\alpha(z)/z\in\CT$ and which we are going to study somewhat closer
in the sequel.

The main result in this note is 

\begin{theorem}
  \label{thm1}
For $f\in\CT$ we have
\begin{equation}
  \label{eq:1}
  \RE\frac{f(\gamma+iy_1)}{f(\gamma+iy_2)}\geq1,\quad
  \gamma\in(-\infty,1],\  0<y_1\leq y_2.
\end{equation}
This relation does not hold, in general, for $\gamma>1$.
\end{theorem}

Theorem~\ref{thm1} has the following immediate consequence.

\begin{corollary}
  \label{cor1}
For $f\in\CT$ and $\gamma\in(-\infty,1]$ fixed, the function $|f(\gamma
+ iy)|$ is monotonically decreasing with $y>0$ increasing.
\end{corollary}

In the case $\gamma=0$ Theorem~\ref{thm1} admits a slight
generalization. It is well-known and easy to verify that $\CT$ is
invariant under the Hadamard product: if
$$
f(z)=\sum_{k=0}^{\infty}a_kz^k,\in\CT,\quad
g(z)=\sum_{k=0}^{\infty}b_kz^k\in\CT,
$$
then also
$$
(f*g)(z):=\sum_{k=0}^{\infty}a_kb_kz^k\in\CT.
$$
\begin{theorem}
  \label{thm2}
For $f,g\in\CT$ we have
$$
\RE\frac{(f*g)(iy)}{f(iy)}\geq1,\quad y>0.
$$
\end{theorem}
And therefore, under the same assumption,  
\begin{equation}
  \label{eq:2}
  |f(iy)|\leq|(f*g)(iy)|,\quad y>0.
\end{equation}
For the polylogarithms and $0<\alpha\leq\beta$
it is clear that
$Li_\beta=Li_\alpha*Li_{\beta-\alpha}$ 
so that we get
\begin{corollary} For $0\leq\alpha<\beta$
$$
|Li_\alpha(iy)|\leq|Li_\beta(iy)|, \quad y>0.
$$
\end{corollary}
This result can also be obtained and even strengthened using  
Corollary~\ref{cor1} and
the deeper relation
$$
\frac{Li_\alpha}{Li_\beta}\in\CT,\quad 0\leq\alpha\leq\beta,
$$
recently established in \cite{rss}.

For a certain subset of $\CT$   
we can go one step beyond Corollary~\ref{cor1},
as far as the behaviour of $|f(iy)|$ for $y>0$ is
concerned.

\begin{theorem}
  \label{thm3}
Let
\begin{equation}
\label{eq:2a}
f(z)=\int_{0}^{1}\frac{\sigma(t)dt}{1-tz},
\end{equation}
where $\sigma\in\CC^1((0,1))$ is positive and with  
$t\sigma'(t)/\sigma(t)$ decreasing.   
Then, for
$w(y):=|f(iy)|$, the function $y w'(y)/w(y)$ decreases with $y>0$ increasing.  
\end{theorem}

Fundamental for the proof of Theorem~\ref{thm3} is the following
result, which is based on a general theorem in \cite{rss}.

\begin{theorem}
  \label{thm4}
Let $f$ be as in Theorem~\ref{thm3}. Then, for $x\in[0,1]$,
$$
\frac{f(z)}{f(xz)}\in\CT.
$$
\end{theorem}

One can show that the conclusion of Theorem~\ref{thm4} is not generally
valid for $f\in \CT$. However, for the functions
$g_\alpha(z):=\frac{1}{z}Li_\alpha(z),\ 
\alpha>0,$ 
we have
$$
g_\alpha(z)=\frac{1}{\Gamma(a)}\int_{0}^{1}\frac{\log^{\alpha-1}(1/t)}{1-tz}dt,
$$
for which the assumptions of Theorem~\ref{thm3} are fulfilled.
Thus both, Theorem~\ref{thm3} and Theorem~\ref{thm4}, apply to $g_\alpha$.  

\section{Proofs} 

We first note that the convex set
$\CT$ satisfies the condition of the main theorem in \cite{r1}, which for
the present case can be stated as follows: 
\begin{lemma}\label{lem1}
Let $\lambda_1, \lambda_2$ be two
continuous linear functionals on $\CT$ and assume that
$0\not\in\lambda_2(\CT)$. Then the range of the functional
$$\lambda(f):=\frac{\lambda_1(f)}{\lambda_2(f)}$$ over $\CT$ equals the
set
$$
\left\{\lambda
\left(
  \frac{\rho}{1-t_1z}+\frac{1-\rho}{1-t_2z}
\right)\ :\ \rho,t_1,t_2\in[0,1]\right\}.
$$
\end{lemma}

\noindent{\bf Proof of Theorem~\ref{thm1}}\ First we note that it is
enough to prove \eqref{eq:1} for $\gamma=1$ only. This is because
$f\in\CT$ implies $f(z-\delta)/f(-\delta)\in\CT$ for all $\delta>0$. 
In Lemma~\ref{lem1} we
choose $\lambda_j(f):=f(1+iy_j),\ j=1,2.$ Since $\IM f(z)>0$ for
$f\in\CT$ and $\IM z>0$, it is clear that
$0\not\in\lambda_2(\CT)$. Lemma~\ref{lem1} now implies that for the
proof of Theorem~\ref{thm1} we only need to show that the expression
$$
\frac{
\dfrac{\rho}{1- t_1-it_1y_1}+\dfrac{1-\rho}{1- t_2-it_2y_1}
}
{\dfrac{\rho}{1- t_1-it_1y_2}+\dfrac{1-\rho}{1- t_2-it_2y_2}},\quad
\rho,t_1,t_2\in[0,1], 
$$
is located in the half-plane $\{w\,:\,\RE w\geq 1\}$. To simplify this
expression we set $\kappa:=(1-\rho)/\rho,\ \tau:=y_1/y_2$. Then our claim is
$$
\RE q(\kappa,y,\tau,t_1,t_2)\geq1,\quad
\kappa\geq0,\ y>0,\ t_1,t_2,\tau\in[0,1],
$$ 
where
$$
q(\kappa,y,\tau,t_1,t_2)= \frac{
\dfrac{1}{1- t_1-i\tau y t_1}+\dfrac{\kappa}{1- t_2-i\tau
  y t_2}
}
{\dfrac{1}{1- t_1-iyt_1}+\dfrac{\kappa}{1- t_2-iyt_2}}.
$$ 
Note that by symmetry we may assume that $t_1\leq t_2$. For fixed
$y,\tau,t_1,t_2$ the values of
$w(\kappa):=q(\kappa,y,\tau,t_1,t_2),\ 
\kappa\geq0$, form a circular arc connecting the points $w(0)=v(t_1)$
and $w(\infty)=v(t_2)$, where
$$
v(t)=\frac{1- t-iyt}{1- t-i\tau yt},
$$
It is easily checked, that under our assumptions for
$y$ and $\tau$  the function $\RE v(t)$ increases with $t\in[0,1]$, 
and, in particular, $\RE v(t)\geq \RE v(0)=1$. This implies
$$
1\leq\RE w(0)\leq\RE w(\infty). 
$$
We will prove that  $\RE w'(0)\geq0$.
Once this done a simple geometric consideration shows that under these
circumstances the circular arc $w(\kappa),\ \kappa\geq0$, cannot leave
the half-plane $\{w\,:\,\RE w\geq1\}$, which then completes the
proof of \eqref{eq:1}.

Calculation yields $$\RE w'(0)=(1-\tau)(t_2-t_1)y^2 \frac ZN$$ where
\begin{eqnarray*}
  Z&=&t_1^*t_2^*(t_2-t_1)+(t_1t_2^*+t_2t_1^*)t_1^*t_2^*\tau +t_1t_2y^2\tau\left(t_1t_2^*+t_2t_1^*-\tau(t_2-t_1)\right),\\
N&=&\left((1-t_1)^2+(t_1y\tau)^2\right)\left((1-t_2)^2+(t_2y\tau)^2\right)\left((1-t_2)^2+(t_2y)^2\right),
\end{eqnarray*}
and $t_j^*:=1-t_j$. Here all terms are non-negative (note that
$$s(\tau):=t_1t_2^*+t_2t_1^*-\tau(t_2-t_1)$$ decreases with $\tau$ and is
therefore not smaller than $s(1)=2t_1t_2^*\geq0$). 
 
It remains to show that \eqref{eq:1} does not hold, in general, for
$\gamma>1$. Let $\gamma=1+\epsilon,\ \epsilon>0$, and choose
$$
f(z):=\frac{1}{1+2\epsilon}+\frac{2\epsilon}{1+2\epsilon}\frac{1}{1-z}\in\CT.
$$
Then, using  $y_1=\epsilon,\ y_2=1,$
$$
\RE\frac{f(\gamma+i\epsilon)}{f(\gamma+i)}=\frac{2\epsilon}{1+\epsilon^2}<1.
$$
\qed

\noindent{}{\bf Proof of Theorem~\ref{thm2}}\ \ If 
$$
g(z)=\int_{0}^{1}\frac{d\mu(t)}{1-tz},
$$
then
$$
\frac{(f*g)(iy)}{f(iy)}=\int_{0}^{1}\frac{f(ity)}{f(iy)}d\mu(t),
$$which is a convex combination of the values of $f(ity)/f(iy)$. By
Theorem~\ref{thm1} these are all in the half-plane $\{w\,:\,\RE
w\geq1\}$.\qed

For the proof of Theorem~\ref{thm4} we need the following result from
\cite{rss}. 

\begin{lemma}
  \label{lem2}
Let $f,g\in\CT$ be represented by
$$
f(z)=\int_0^1\frac{\phi(t)dt}{1-tz},\quad g(z)=\int_0^1\frac{\psi(t)dt}{1-tz}
$$
with non-negative Borel functions $\phi,\psi$ on $(0,1)$. If 
$\phi(t)\psi(s)\geq\phi(s)\psi(t)$ holds for all $0<s<t<1$, then $f/g\in\CT$.
\end{lemma}

\noindent{\bf Proof of Theorem~\ref{thm4}}\ We have 
$$
f(xz)=\int_{0}^{1}\frac{\sigma(t)dt}{1-txz}
=\int_{0}^{1}\frac{\sigma^*(t)dt}{1-tz}, 
$$
with
$$
\sigma^*(t):=\left\{
  \begin{array}{ll} \displaystyle 
    \frac{1}{x}\sigma(t/x),&0<t\leq x,\\[4mm]
0,&x<t<1.
  \end{array}
\right.
$$
The condition 
\begin{equation}
\label{eq:3}
\sigma(t)\sigma^*(s) \geq \sigma(s)\sigma^*(t),\quad 0<s<t<1,
\end{equation}
is immediately fulfilled if $t>x$. Otherwise we are left with
$$
\sigma(t)\sigma(s/x) \geq \sigma(s)\sigma(t/x),\quad 0<s<t\leq x.
$$
This requires that $\sigma(t)/\sigma(t/x)$ increases with $t$. Taking
logarithms and differentiating w.r.t. the variable $t$, we find as a
necessary and sufficient condition for \eqref{eq:3} that
$t\sigma'(t)/\sigma(t)$ decreases for $t$ increasing. The result
follows now from Lemma~\ref{lem2}.\qed
\vspace{2mm}

\noindent{\bf Proof of Theorem~\ref{thm3}}\ We apply
Theorem~\ref{thm1} to the function $F$ of Theorem~\ref{thm4}. 
Then, for $x,\tau\in(0,1)$, we get
$$
\left|
  \frac{f(iy\tau)f(iyx)}{f(iyx\tau)f(iy)}
\right|\geq1,\quad y>0.
$$
Taking logarithms we obtain
$$
(\log w(y)-\log w(x y))-(\log w(\tau y)-\log w(x\tau y))\leq0,
$$
Dividing by $1-x$ and letting $x\rightarrow 1-0$ yields
$$
\frac{y w'(y)}{w(y)}\leq\frac{\tau y w'(\tau y)}{w(\tau y)},
$$
which implies the assertion.\qed

\end{document}